\def\VersionName{\jobname.tex}

\documentclass[12pt]{article}
\usepackage{amsmath,amssymb,amsthm}

\normalsize

\date{}

\numberwithin{equation}{section}
\newtheorem{theorem}{Theorem}[section]
\newtheorem{corollary}{Corollary}[section]
\newtheorem{lemma}{Lemma}[section]
\newtheorem{remark}{Remark}[section]

\def\reals{\mathbb{R}}

\def\angles{\reals_{2\pi}}

\def\x{\boldsymbol{x}}
\def\0{{\mathbf 0}}
\def\1{{\mathbf 1}}
\def\btheta{\boldsymbol{\theta}}

\def\S{\mathcal{I}_k}
\def\R{\mathcal{R}}
\def\Rc{{\mathcal{C}}}
\def\RRc{{\R{+}\Rc}}
\def\radius{D}
\def\Group{P}
\def\Re{\mathrm{Re}}

\def\({\bigl(}
\def\){\bigr)}
\def\midcolon{\mathrel{:}}
\def\abs#1{\mathopen|#1\mathclose|}

\def\dfrac#1#2{\lower0.15ex\hbox{\large$\frac{#1}{#2}$}}

\title{Asymptotic enumeration of correlation-immune boolean functions}

\author{
E. Rodney Canfield\thanks
 {Research supported by the NSA Mathematical Sciences Program} \\
\small Department of Computer Science\\[-0.8ex]
\small University of Georgia\\[-0.8ex]
\small Athens, GA 30602, USA\\[-0.3ex]
\small\texttt{erc@cs.uga.edu}
\and
Zhicheng Gao\\
\small School of Mathematics and Statistics\\[-0.8ex]
\small Carleton University, Ottawa, Canada\\[-0.3ex]
\small\texttt{zgao@math.carleton.ca}
\and
Catherine Greenhill\\
\small School of Mathematics and Statistics\\[-0.8ex]
\small University of New South Wales\\[-0.8ex]
\small Sydney, Australia 2052\\[-0.3ex]
\small\texttt{csg@unsw.edu.au}
\and
Brendan D. McKay\thanks{Research supported by the
 Australian Research Council}\\
\small School of Computer Science\\[-0.8ex]
\small Australian National University \\[-0.8ex]
\small Canberra, ACT Australia\\[-0.3ex]
\small\texttt{bdm@cs.anu.edu.au}
\and
Robert W. Robinson \\
\small Department of Computer Science\\[-0.8ex]
\small University of Georgia\\[-0.8ex]
\small Athens, GA 30602, USA\\[-0.3ex]
\small\texttt{rwr@cs.uga.edu}
}

\begin{document}

\setlength{\abovedisplayskip}{0.7\abovedisplayskip}
\setlength{\belowdisplayskip}{0.7\belowdisplayskip}
\setlength{\abovedisplayshortskip}{0.8\abovedisplayshortskip}
\setlength{\belowdisplayshortskip}{0.8\belowdisplayshortskip}

\maketitle

\newpage

\begin{abstract}
A boolean function of $n$
boolean variables is \textit{correlation-immune\/}
of order~$k$ if the function value is uncorrelated with the values of
any $k$ of the arguments.  Such functions are of considerable
interest due to their cryptographic properties, and are also
related to the orthogonal arrays of statistics and the balanced
hypercube colourings of combinatorics.
The \textit{weight\/} of a boolean function
is the number of argument values that produce a function value
of~1.  If this is exactly half the argument values, that is,
$2^{n-1}$ values, a
correlation-immune function is called \textit{resilient}.

An asymptotic estimate of the number $N(n,k)$ of $n$-variable
correlation-immune
boolean functions of order~$k$ was obtained in 1992 by Denisov
for constant~$k$.  Denisov repudiated that estimate in
2000, but we will show that the repudiation was a mistake.

The main contribution of this paper is an asymptotic estimate
of $N(n,k)$ which holds if $k$ increases with $n$ within
generous limits and specialises to functions with a given weight,
including the resilient functions.
In the case of $k=1$, our estimates are valid for all weights.
\end{abstract}

\section{Introduction}\label{section:Intro}

Let $n,k,q$ be integers satisfying $1\le k\le n$ and $0\le 2^kq\le 2^n$,
and define $\lambda=2^kq/2^n$.\linebreak
A \textit{correlation-immune boolean function} of $n$ variables,
order $k$ and weight~$2^kq$ is a boolean-valued function
of $n$ boolean variables with this property:
if any $k$ arguments are given arbitrary values, exactly the
fraction $\lambda$ of the $2^{n-k}$ possible assignments to
the remaining arguments give a function value
of~$1$.  (See for example~\cite{Maitra, Roy, Schneider} 
and~\cite[Chapter 4]{CS}.)
Let $N(n,k,q)$ denote the number of such functions.
An important special case is the \textit{resilient functions\/},
which have $\lambda=\tfrac12$.
Correlation-immune functions, and in particular the resilient functions,
have desirable cryptographic properties:  see for 
example~\cite{Carlet-chapter, Roy}.
In this paper we will derive an asymptotic estimate of $N(n,k,q)$
for a wide range of $k$ and $q$ values, and deduce an asymptotic formula
for the sum  $N(n,k) = \sum_q N(n,k,q)$, which is the number
of correlation-immune boolean functions of $n$ variables and
order $k$.

An $n$-variable boolean function can be represented as a matrix of
$n$ columns over $\{0,1\}$ whose rows consist of those argument lists
which give the function value~1.   A correlation-immune boolean function
of $n$ variables, order $k$ and weight $2^k q$  gives rise to
a matrix with $2^k q$ distinct rows and $n$ columns, such that 
in any set of $k$ columns each of the $2^k$ possible 0-1 patterns appears
exactly $q$ times.
In statistics, such a matrix is called an \textit{orthogonal array}
of 2~levels, $n$~variables, $2^kq$ runs, and strength~$k$;
see~\cite{Hedayat} for a detailed exposition.
Since the $2^kq$ rows are by definition distinct, and
permuting the rows does not change the associated function, there is an
uninteresting ratio of $(2^kq)!$ between the number of matrices and
the number of functions. We will work with functions rather than
matrices.

The special case $k=1$ has also been studied under the name
of \textit{balanced colourings of a hypercube}.  These are
placements of equal weights on some of the vertices of a
hypercube such that the centroid is at the center of the
hypercube.  Exact enumerations have been found
in this case~\cite{Palmer,Zhang}, but they do not appear
suitable for asymptotics.

Early papers on the number of correlation-immune functions focussed on
the case $k=1$.   Upper and lower bounds for $N(n,1)$ were
given in~\cite{Maitra, mitchell, PLSK, YG} but these do not
appear as sharp as the bounds given by Bach~\cite{Bach}.

The case of general $k$ was first considered by
Schneider~\cite{Schneider}, who gave upper bounds for
$N(n,k,q)$ as well as for $N(n,k)$.
For large $k$ an improved
upper bound is given by Carlet and Klapper~\cite{CK},
both for $N(n,k)$ and for the resilient functions of order~$k$.
Carlet and Gouget~\cite{Carlet} gave an upper bound 
for the number of resilient functions of order~$k$, which
improves upon Schneider's bound for $k> n/2-1$ and which partially
improves upon the upper bound of~\cite{CK}.
Tarannikov~\cite{tarannikov} proved that when $c$ is a fixed positive
integer, the function $N(n,n-c)$ is bounded above by a polynomial in $n$.
Exact expressions for $N(n,n-c)$ when $c=1,2,3$ are also given
in~\cite[Theorem 3]{tarannikov}.
(See also~\cite{TK}.)

The first asymptotic enumeration of correlation-immune functions
was achieved by Denisov.  
Define
\[
      M = \sum_{j=0}^k \binom{n}{j}\mbox{~~and~~}
      Q = \sum_{j=1}^k j\mkern1mu\binom{n}{j}.
\]

\begin{theorem}[Denisov \cite{Denisov1}]\label{denisov1}
If $k\ge 1$ is a constant integer then, as $n\to\infty$,
\[
    N(n,k)
      \sim
    2^{2^n+Q-k}(2^{n-1}\pi)^{-(M-1)/2}.\quad\qedsymbol
\]
\end{theorem}
Denisov's formula for $N(n,1)$ was refined by Bach~\cite{Bach},
who showed that an asymptotic expansion for $N(n,1)$ exists and 
calculated the first few terms of it.  

In a later paper~\cite{Denisov2}, Denisov repudiated
Theorem~\ref{denisov1} and proposed a different value.
However, we will show that Denisov's repudiation was a
mistake, and Theorem~\ref{denisov1} is correct.  More
discussion of \cite{Denisov2} is given in Section~\ref{finalsection}.

\bigskip

We now state our results. Define
\[ A = \lambda(1-\lambda).\]
In addition to common asymptotic notations like $O(\cdot)$,
we use $\omega(f(n))$ to represent any function $g(n)$ such
that $g(n)/f(n)\to\infty$ as $n\to\infty$.

\begin{theorem}\label{main}
Consider a sequence of triples $(n,k,q)$ of positive integers
such that $n\to\infty$ and
\begin{equation}\label{kbounds}
\omega\(2^{5k}n^{6k+3}M^3\)\le q\le 2^{n-k}-\omega\(2^{5k}n^{6k+3}M^3\).
\end{equation}
Then
\begin{equation}\label{theanswer}
N(n,k,q) 
 = 
2^Q\,
            \(\lambda^{\lambda }(1-\lambda)^{1-\lambda}\)^{-2^n}\,
\(\pi A \, 2^{n+1}\)^{-M/2}\,
\(1+O(\eta(n,k,q))\),
\end{equation}
where
$\eta(n,k,q) =
  2^{-n/2+3k}n^{3k+3/2}M^{3/2}\lambda^{-1/2}(1{-}\lambda)^{-1/2}
  = o(1)$.
\end{theorem}

\begin{remark}\label{remark1}
Given a function $g$ in the class counted by $N(n,k,q)$, we
can form another, namely $1-g$, counted by $N(n,k,2^{n-k}{-}q)$.
This complementation operation is a bijection which
exchanges $q$ with $2^{n-k}{-}q$ and $\lambda$ with $1{-}\lambda$.
This means, for example, that we can assume $\lambda\le\tfrac12$
in our proof when it is convenient.
\end{remark}

\begin{remark}\label{remark2}
By Stirling's formula, $\log M = o(n)$ whenever $k=O(n/\log n)$.
{}From this it follows that~\eqref{kbounds} is non-vacuous whenever
\begin{equation}\label{krange}
1\le k \le \biggl(\frac{\log 2}{6}-\varepsilon\biggr)\frac{n}{\log n}
\end{equation}
for some $\varepsilon>0$.
\end{remark}

\begin{corollary}\label{resilient}
  If\/ $k=k(n)$ satisfies~\eqref{krange}
then, as $n\to\infty$, the number of $k$-resilient boolean
functions of $n$ variables~is
\[
   2^{2^n+Q}(2^{n-1}\pi)^{-M/2} \( 1+O(2^{-n/2+3k}n^{3k+3/2}M^{3/2})\)
   \sim 2^{2^n+Q}(2^{n-1}\pi)^{-M/2}.
\]
 
\end{corollary}

\begin{corollary}\label{rightfirsttime}
  If\/ $k=k(n)$ satisfies~\eqref{krange}
then, as $n\to\infty$, the number of order $k$ correlation-immune
boolean functions of $n$ variables is
\begin{align}
    N(n,k) &= 2^{2^n+Q-k}(2^{n-1}\pi)^{-(M-1)/2}
    \( 1+O(2^{-n/2+3k}n^{3k+3/2}M^{3/2})\)\label{Nsum}\\
    &\sim 2^{2^n+Q-k}(2^{n-1}\pi)^{-(M-1)/2}.\notag
\end{align}
\end{corollary}

Corollary~\ref{rightfirsttime} shows that Denisov's result
Theorem~\ref{denisov1} is true, despite his later retraction.

In Section~\ref{integralsection}, we write $N(n,k,q)$ as an
integral in many complex dimensions.
In Section~\ref{domainsection} we identify the points where the
integrand has maximum magnitude and define a region
$\RRc$ consisting of a small hypercuboid surrounding each of
those points.
The integral is then bounded outside $\RRc$ in
Section~\ref{outofboxsection}
and estimated inside $\RRc$ in Section~\ref{inboxsection}.
The proof of Theorem~\ref{main} is completed in
Section~\ref{proofsection} where we also prove
Corollaries~\ref{resilient} and~\ref{rightfirsttime}.  In the final sections
we consider some additional topics including a closer look
at the case $k=1$ and a connection with Hadamard matrices.

\section{The desired quantity
              as a complex integral}\label{integralsection}

Define $[n]=\{1,2,\dots,n\}$ and
$\S=\{S\in 2^{[n]}\midcolon |S|\le k\}$.
We will identify $N(n,k,q)$ as the constant term in a generating
function over the $M$ variables $\{x_S \midcolon S\in\S\}$.
Let $\x$ denote a vector of all these variables,
in arbitrary order.
For $\radius=\lambda/(1-\lambda)$,
define the rational function $F(\x)$ by
$$
F(\x)=\prod_{\alpha\in\{\pm 1\}^n}\biggl(1+\radius\prod_{S\in\S}
\! x_{S}^{\alpha_S} \biggr),
$$
where
$$
\alpha_S=\prod_{j\in S}\,\alpha_j
$$
for each $S$ (including the case $\alpha_\emptyset=1$).
The value of $\radius$ is determined by a saddle point condition, as
will become apparent in Section~\ref{inboxsection}.

\begin{lemma}\label{Nlemma}
$N(n,k,q)$ is the constant term of\/ $(\radius x_\emptyset)^{-2^kq}F(\x)$.
\end{lemma}

\begin{proof}
For a boolean function $g(y_1,\ldots, y_n)$, the
\emph{Walsh transform} of $g$ is the real-valued function
$\hat{g}$ over $\{ 0,1\}^n$ defined by
\[ \hat{g}(w_1,\ldots, w_n) = \sum_{(y_1,\ldots, y_n)\in\{ 0,1\}^n}
                 g(y_1,\ldots, y_n)\, (-1)^{w_1 y_1 + \cdots + w_n y_n}.
                 \]
Given $\alpha\in\{ \pm 1\}^n$, form $\bar{\alpha}\in\{ 0,1\}^n$
from $\alpha$ by changing each 1 entry into 0 and each $-1$
entry into 1.

For $S\in\S$, let $w_S\in\{ 0,1\}^n$ be the
characteristic vector of $S$.
Then, given a vector $\alpha\in\{\pm 1\}^n$ and any $S\in\S$, we
have
\[ \alpha_S = (-1)^{\bar{\alpha}\cdot w_S}.\]
We can view $F(\x)$ as the sum of $2^{2^n}$ terms,
with one term for each boolean function $g$ of $n$ variables.
Specifically, the term corresponding to a boolean function
$g:\{0,1\}^n\to\{0,1\}$ is exactly
\[ \prod_{\stackrel{\alpha\in\{ \pm 1\}^n}{ g(\bar{\alpha})=1}}
      \biggl(
                  \radius \prod_{S\in \S} x_S^{\alpha_S}\biggr)
     = \radius^{\hat{g}(w_\emptyset)}
             \prod_{S\in\S} x_S^{\hat{g}(w_S)}.\]
By the spectral characterisation of correlation-immune
functions~\cite{Sarkar, Xiao}, the boolean function $g$ is
correlation-immune of order $k$ if and only if
$\hat{g}(w_S) = 0$ for all $S\in\S \setminus \{ \emptyset\}$.
Moreover, the functions counted by $N(n,k,q)$ have
$\hat{g}(w_\emptyset) = 2^k q$.  Therefore
the coefficient of the monomial $x_\emptyset^{2^k q}$ in $F(\x)$
is exactly $\radius^{2^k q}\, N(n,k,q)$.
\end{proof}

By Cauchy's integral formula, it follows from Lemma~\ref{Nlemma} that
$$
N(n,k,q)
=\frac{1}{(2\pi i)^M \radius^{2^k q}}
\oint \cdots \oint \frac{F(\x)}
         {x_\emptyset^{2^k q}\prod_{S\in\S}x_S} \,d\x,
$$
where each $x_S$ is integrated anticlockwise around a circle of radius~1
centred at the origin.
Now introduce variables $\theta_S$ ($S\in\S$)
and the $M$-dimensional vector $\btheta$
of the $\theta_S$ variables in arbitrary order.
Change variables from $\x$ to $\btheta$ using
$x_S=e^{i\theta_S}$ for each~$S$.
Then
\begin{align}
N(n,k,q)
& =\frac{(1+\radius)^{2^n}}{(2\pi)^{M}\radius^{2^kq}}\, I(n,k,q), 
   \label{Idef}\\
\intertext{where}
I(n,k,q) &=
\int_{-\pi}^{\pi} \cdots \int_{-\pi}^{\pi} G(\btheta) \,d\btheta, \notag\\
G(\btheta) &= e^{-i2^kq\theta_\emptyset}\!\!\prod_{\alpha\in\{\pm1\}^n}
  \!\frac{1+\radius e^{if_{\alpha}(\btheta)}}{1+\radius}, \label{Gdef}\\
\intertext{and}
f_{\alpha}(\btheta) 
 &= \sum_{S\in\S} \alpha_S \theta_S.
\end{align}

The elements of $\btheta$ belong to the set $\angles$ of real numbers
modulo $2\pi$.
In this set, addition, and multiplication by integers,
have their usual meanings.
We use $\equiv$ to indicate equality in $\angles$.  For
example, $\theta \equiv 0$ means that $\theta$ is the element of
$\angles$ corresponding to the real number $2\pi t$ for any
integer $t$.  Also let 
\[ z:\angles\to(-\pi,\pi]\]
be the standard mapping of
$\angles$ onto the real interval $(-\pi,\pi]$
and define the absolute value
$d(\theta)=\abs{z(\theta)}$ for any $\theta\in\angles$.
Clearly $d(\cdot)$ satisfies the triangle inequality:
$d(\theta+\theta')\le d(\theta)+d(\theta')$.

\section{Analysis of the domain of integration}\label{domainsection}

The integrand $G(\btheta)$ defined in~\eqref{Gdef} has modulus
at most~1. We will later show that the value of the integral
$I(n,k,q)$ comes mostly from the near vicinity of those points
where equality occurs, so our next task will be to identify those
points.  Define
\[
 \Rc = \bigl\{\, \btheta\in\angles^M \, \midcolon \, \abs{G(\btheta)}=1\,\bigr\}.
\]

\begin{lemma}\label{Rclemma}
\begin{equation}\label{Rcdef}
   \Rc = \Bigl\{\, \btheta\in\angles^M \, \midcolon \,
      2^{|S|}\mkern-8mu\sum_{T\in\S, T\supseteq S}\!\!\!\theta_T
         \equiv 0\,\mbox{ for each $S\in\S$}\Bigr\},
\end{equation}
and moreover $|\Rc|=2^Q$.
\end{lemma}

\begin{proof}
Throughout the proof we work in $\angles$.
For $1\le j\le n$, define the linear difference operator $\delta_j$ by
\[
   \delta_j f_{(\alpha_1,\ldots,\alpha_j,\ldots,\alpha_n)}
   = f_{(\alpha_1,\ldots,\alpha_j,\ldots,\alpha_n)}
  - f_{(\alpha_1,\ldots,\alpha_{j-1},-\alpha_j,\alpha_{j+1},\ldots,\alpha_n)}.
\]
For $S\in\S$, define
$\delta_S=\prod_{j\in S}\delta_j$, noting that the product is
commutative.
{}From the definition of $f_\alpha(\btheta)$ we can easily prove by
induction on $|S|$ that
\begin{equation}\label{Sdiff}
\delta_S f_{\alpha}(\btheta)
= 2^{|S|}\mkern-8mu \sum_{T\in\S, T\supseteq S} \alpha_T \theta_T.
\end{equation}

Since 
\[ \biggl| \frac{1+De^{ix}}{1+D}\biggr| = \frac{\sqrt{1+2D\cos(x) + D^2}}{1+D} \leq 1,\]
a necessary and sufficient condition for $\btheta\in\Rc$
is that $f_\alpha(\btheta)\equiv 0$ for all~$\alpha\in \{\pm 1\}^n$.

Suppose that $\btheta\in\Rc$.
Then, since $f_\alpha(\btheta)\equiv 0$ for all~$\alpha$,
the difference $\delta_S f_{\alpha_0}(\btheta)$
satisfies $\delta_S f_{\alpha_0}(\btheta) \equiv 0$
for all $S\in\S$, where $\alpha_0 = (1,1,\ldots, 1)$.
By~\eqref{Sdiff} we conclude that $\btheta$ lies in the
set $\Rc^*$ given by the right hand side of~\eqref{Rcdef},
and hence $\Rc\subseteq \Rc^*$.
Conversely, if $\btheta\in\Rc^*$ then every $f_\alpha(\btheta)\equiv 0$
since
\[
   f_\alpha(\btheta)
   \equiv \Bigl({\textstyle\prod_{\{j\midcolon\alpha_j=-1\}} (1-\delta_j)}\Bigr)
    f_{\alpha_0}(\btheta).
\]
Therefore, $\Rc=\Rc^*$.

Since the set of equations in~\eqref{Rcdef} is triangular, we
can find all solutions by choosing each $\theta_S$ in order of
decreasing~$|S|$.  
There are exactly
$2^{|S|}$ choices for $\theta_S$, so the total number of solutions is
$|\Rc|=2^Q$.
\end{proof}

As noted in Remark~\ref{remark1},
we will assume that $\lambda\le\tfrac12$ without losing generality.
Let $\Delta$ be the positive number defined by
\[
\Delta = 2^{-n/2+k+3} \lambda^{-1/2} n^{k+1/2} M^{1/2}.
\]
The left side of~\eqref{kbounds} is equivalent to
\begin{equation}\label{Deltalimit}
  \Delta = o\(2^{-2k}n^{-2k-1}M^{-1}\).
\end{equation}
Let $\R$ be the subset of $\angles^M$ defined by 
\[
 \R = \bigl\{ \,\btheta\in \angles^M \, \midcolon  \,
 d(\theta_S)\le\Delta(2n)^{-|S|} \,\, \mbox{ for all } \,\, S\in\S
 \, \bigr\}.
\]
This is a hypercuboid centred at the origin.
Denote the union of $2^Q$ copies of~$\R$ centred at
the points in $\Rc$ by 
$$\RRc= \bigcup_{\btheta^*\in\Rc}\{\R+\btheta^*\}\subseteq \angles^M.$$

Since all the elements of vectors in $\Rc$ are integer multiples of
$2\pi/2^k$, it follows from~\eqref{Deltalimit} that these
copies are disjoint. The region $\RRc$ includes all the points
where $\abs{G(\btheta)}$ is maximal; we will prove in the following
sections that in fact
it includes all the points which contribute substantially to
$I(n,k,q)$.  

\section{The integral outside the critical region}
\label{outofboxsection}

\begin{lemma}\label{outofbox}
If\/ the conditions of Theorem~\ref{main} are satisfied
and $\lambda \leq \tfrac12$ then
\[
  \int_{(\RRc)^c} \abs{G(\btheta)}\,d\btheta 
 < (2\pi)^M \exp\(-\tfrac{4}{5}nM\),
\]
where $(\RRc)^c = \angles^M \setminus (\RRc)$.
\end{lemma}

\begin{proof}
Fix $\btheta\in (\RRc)^c$.
First we show that there exists some set $S_0 = S_0(\btheta)\in\S$ such that
\begin{equation}
 d(\delta_{S_0}f_{\alpha}(\btheta))>(2-e^{1/2})\Delta n^{-|S_0|}
\label{step1}
\end{equation}
for all $\alpha\in\{\pm1\}^n$.
Define
$\btheta^*=\btheta^*(\btheta)\in\Rc$ recursively, as follows:  
starting with sets $S\in\S$ with $|S|=k$, and then proceeding to smaller $k$, 
choose $\theta^*_S\in\angles$ such that
$2^{|S|}\sum_{T\in\S, T\supseteq S} \theta^*_T \equiv 0$ 
and $d(\theta_S-\theta^*_S)$ is minimal over
all such choices of $\theta^*_S$.
(Break ties arbitrarily.)
Since $\btheta\notin\RRc$, there is a set $S_0\in\S$ of maximum
cardinality such that
\begin{equation}\label{wantit}
d(\theta_{S_0}-\theta^*_{S_0}) > (2n)^{-|S_0|}\Delta.
\end{equation}
By the maximality of $S_0$ we have
\begin{align}\label{tailbound}
\sum_{T\in \S, T\supset S_0} \!\! d(\theta_T-\theta^*_T)
 &\le \sum_{j\ge 1} \binom{n}{j} \Delta(2n)^{-|S_0|-j} \notag\\
 &\le \Delta (2n)^{-|S_0|}\sum_{j\ge 1} \frac{2^{-j}}{j!} \notag\\
 & = (e^{1/2}-1)\Delta(2n)^{-|S_0|}.
\end{align}
Now take any $\alpha\in\{\pm1\}^n$ and write, using (\ref{Sdiff}),
\[
 \delta_{S_0} f_\alpha(\btheta) \equiv 2^{|S_0|} 
  \negthickspace\sum_{T\in\S, T\supseteq S_0} \negthickspace
  \alpha_T\theta_T
  \equiv \Sigma_1 + \Sigma_2 + U
  \]
  where
  \begin{align*} 
 \Sigma_1 &\equiv 
       2^{|S_0|} \negthickspace\sum_{T\in \S, T\supseteq S_0}
         \negthickspace\alpha_T\theta^*_T,\\
 \Sigma_2 &\equiv 
       2^{|S_0|} \negthickspace\sum_{T\in \S, T\supset S_0}
         \negthickspace\alpha_T (\theta_T-\theta^*_T),\\
       U &\equiv 
       2^{|S_0|} \alpha_{S_0}(\theta_{S_0}-\theta^*_{S_0}).
\end{align*}
Since $\btheta^*\in\Rc$, \eqref{Rcdef} implies that $\Sigma_1\equiv 0$.
Next, since $d(\alpha_T\theta_T)=d(\theta_T)$,
\eqref{tailbound} implies that
\[ d(\Sigma_2)\leq (e^{1/2}-1)\Delta n^{-|S_0|}.\]
Finally, 
\[ d(U)>\Delta n^{-|S_0|}, \]
by~\eqref{wantit} and the fact that 
$d(\theta_{S_0}-\theta_{S_0}^*) < 2^{-|S_0|}\pi$.
Therefore, using the triangle inequality,
\[
  d(\delta_{S_0}f_{\alpha}(\btheta))
  = d(U + \Sigma_2)
  \geq d(U) - d(\Sigma_2)
   > (2-e^{1/2})\Delta n^{-|S_0|}.
\]
Since $\alpha\in\{ \pm 1\}^n$ was arbitrary, this 
establishes the existence of the desired set $S_0$.

Next, partition the set $\{ \pm 1\}^n$ into $2^{n-|S_0|}$
parts, each of size $2^{|S_0|}$, such that 
two vectors $\alpha, \alpha'$ belong to the same part if
and only if they agree in every coordinate $j\not\in S_0$.
Let $\Group$ be an arbitrary part of the partition.  
For any $\alpha\in\Group$,
the difference $\delta_{S_0}f_{\alpha}(\btheta)$ is a linear combination, with
coefficients $\pm 1$, of the elements of the set
$\{ f_{\alpha'}(\btheta) \midcolon  \alpha'\in \Group\}$.
Therefore, by~\eqref{step1} and using the triangle inequality,
\begin{equation}
\label{step2}
(2-e^{1/2})\Delta n^{-|S_0|} < d(\delta_{S_0}f_{\alpha}(\btheta))
  \le \sum_{\alpha'\in\Group} d(f_{\alpha'}(\btheta)).
\end{equation}
As $1 - \cos x \leq 2x^2/\pi^2$ for $-\pi \leq x\leq \pi$,
we find that for all $x\in\mathbb{R}$,
\begin{align*}
 \left| \frac{1+\radius e^{ix}}{1+\radius}
\right|^{\,2} &= 1 - \frac{2\radius(1-\cos x)}{(1+\radius)^2} \\
 &\le
\exp\Bigl(-\frac{4\radius \, d(x)^2}{(1+\radius)^2 \pi^2} \Bigr) \\
&= \exp\Bigl(-\frac{4\lambda(1-\lambda)}{\pi^2}\, d(x)^2\Bigr)\\
&\leq \exp\Bigl(-\frac{2\lambda}{\pi^2}\, d(x)^2\Bigr)
\end{align*}
using the assumption $\lambda \le \tfrac12$ for the last inequality.
Thus, using the Cauchy-Schwarz inequality and~\eqref{step2},
\begin{align*}
\prod_{\alpha\in\Group}\,
  \biggl| \frac{1+\radius e^{if_{\alpha}(\btheta)}}{1+\radius} \biggr| &\le
  \exp\biggl(-\frac{\lambda}{\pi^2}
    \sum_{\alpha\in\Group}d(f_{\alpha}(\btheta))^2\biggr) \\
&\le \exp\biggl(-\frac{\lambda}{\pi^2\abs\Group}
  \biggl(\, \sum_{\alpha\in\Group}
     d(f_{\alpha}(\btheta))\biggr)^{\!\!2}\,\biggr)\\
&\le \exp\biggl(-\frac{\lambda}{\pi^2}\, 2^{-|S_0|}\, 
                     ((2-e^{1/2})\Delta n^{-|S_0|})^2\biggr).
\end{align*}
Since there are $2^{n-|S_0|}$ parts in the partition,
taking the product over all parts and applying the
definition of $\Delta$ gives
\begin{align*}
\abs{G(\btheta)}
&\le \exp\(
-(2-e^{1/2})^2 \pi^{-2}
2^{2k-2|S_0|+6} n^{2k-2|S_0|+1}M \)\\
&\le \exp\( -2^6\, (2-e^{1/2})^2 \pi^{-2} nM\),
\end{align*}
as $|S_0|\leq k$.
Finally we note that $2^6\, (2-e^{1/2})^2\, \pi^{-2}>\tfrac 45$,
so we have
\[
   \abs{G(\btheta)} < \exp\(-\tfrac45 nM\).
\]
As this inequality holds for any $\btheta\notin\RRc$  
and the volume of $(\RRc)^c$ is at most $(2\pi)^M$, 
the proof is complete.
\end{proof}

\section{The integral inside the critical region}
\label{inboxsection}

\begin{lemma}\label{inthebox}
If\/ the conditions of Theorem~\ref{main} are satisfied and
$\lambda\le\tfrac12$ then
\[
\int_\R G(\btheta)\,d\btheta =
  \biggl( \frac{2\pi}{\lambda(1-\lambda)2^n}
           \biggr)^{\!\!M/2}
      \(1 + O(2^{5k/2}n^{3k+3/2}M^{3/2}q^{-1/2})\).
\]
\end{lemma}

\begin{proof}
Let $\btheta=(\theta_S)_{S\in\S}\in\R$.
In this section we perform expansions that are valid in $\reals$
rather than $\angles$, so we identify $\btheta$ with
$\(z(\theta_S)\)_{S\in\S}$.
Since
\[
   \exp\biggl( i\sum_{S\in\S} \alpha_S\theta_S \biggr)
   = 
   \exp\biggl( i\sum_{S\in\S} \alpha_S z(\theta_S) \biggr),
\]
this identification has no effect on $G(\btheta)$.
Also note that
\begin{equation}\label{falphabound}
  \abs{f_{\alpha}(\btheta)} = 
  \biggl|\sum_{S\in\S} \alpha_S z(\theta_S)\biggr|
  \le \Delta \sum_{j=0}^k \binom{n}{j}(2n)^{-j} \leq  e^{1/2}\Delta.
\end{equation}
Define
\[ h(x) = \log\biggl(\frac{1+\radius  e^{ix}}{1+\radius }\biggr).\]
By Taylor's Theorem with the integral form of the remainder
(which also holds for complex-valued functions),
\[
h(f_\alpha(\btheta))
= i\,\frac{\radius }{1+\radius } \, f_{\alpha}(\btheta)
        - \dfrac12\frac{\radius }{(1+\radius )^2}\, f_{\alpha}(\btheta)^2
  + R(f_\alpha(\btheta)) 
\]
where  
\begin{equation}
\label{remainder}
 R(f_\alpha(\btheta)) = \int_{0}^{f_\alpha(\btheta)} \dfrac{1}{2} 
            h'''(t)(f_\alpha(\btheta)-t)^2 dt.
\end{equation}
Now $\cos(\cdot)$ is unimodal on $[-e^{1/2}\Delta,e^{1/2}\Delta]$
by (\ref{Deltalimit}).
Therefore for $|t|\leq e^{1/2}\Delta$ we have
\[
     |h'''(t) | 
      = \frac{D\sqrt{1-2D\cos(t) + D^2}}{(1+2D\cos(t)+D^2)^{3/2}}
   \leq \, D \leq 2\lambda
\]
using the assumption that $\lambda \leq \tfrac12$.
Hence by~\eqref{falphabound} and~\eqref{remainder}, 
\begin{equation*}
 \left|R(f_\alpha(\btheta))\right| 
                     \leq
    \frac{\lambda\, e^{3/2} \Delta^3 } 
              {3}
  \leq 2 \lambda \Delta^3. 
\end{equation*} 
Then
\begin{align*}
 G(\btheta)
&= \exp\biggl(-i2^kq\theta_\emptyset + 
\sum_{\alpha\in\{ \pm 1\}^n} 
   \left( i\frac{\radius }{1+\radius }
               f_{\alpha}(\btheta)
   - \dfrac12\frac{\radius }{(1+\radius )^2}
                 f_{\alpha}(\btheta)^2 + R(f_\alpha(\btheta)) \right)\biggr)  \\
 &= \exp\biggl(\,\,
 \sum_{\alpha\in\{ \pm 1\}^n}
   \left( - \dfrac12\frac{\radius }{(1+\radius )^2} f_{\alpha}(\btheta)^2
    + R(f_\alpha(\btheta)) \biggr)\right)  \\
   &= \exp(a(\btheta)) \,
   \exp\biggl( -\dfrac{1}{2}\, A 
                \sum_{\alpha\in\{ \pm 1\}^n} f_{\alpha}(\btheta)^2
   \biggr)
\end{align*}
where
\[ a(\btheta) = \sum_{\alpha\in \{ \pm 1\}^n } R(f_\alpha(\btheta)).\]
The vanishing of the linear terms explains our choice of~$\radius $.
Note that $a(\btheta)$ is a complex number which is bounded in
modulus by
\begin{equation}
\label{a-bounds}
 |a(\btheta)| \leq \lambda\, 2^{n+1} \Delta^3. 
\end{equation}

Next, note that the reflection $\btheta\mapsto -\btheta$ preserves the
region $\R$ and maps $G(\btheta)$ to its complex conjugate.  
It follows that
\[ \int_{\R} G(\btheta)\, d\btheta\]
is real, and therefore is equal to the integral of the real part of
its integrand.   Hence
\begin{align*}
\int_{\R} G(\btheta) \, d\btheta &= 
     \int_{\R} \Re(\exp(a(\btheta)))\, \exp \biggl( -\dfrac{1}{2}\, A 
                \sum_{\alpha\in\{ \pm 1\}^n} f_{\alpha}(\btheta)^2 
                                      \biggr) \, d\btheta  \notag \\
        &= \Re(\exp(a(\btheta_0)))\, \int_{\R} \exp \biggl( -\dfrac{1}{2}\, A 
                \sum_{\alpha\in\{ \pm 1\}^n} f_{\alpha}(\btheta)^2
                                       \biggr) \, d\btheta 
\end{align*}
for some $\btheta_0\in\R$, using the Intermediate Value Theorem.

Since $\lambda 2^n \Delta^3 = o(1)$ using (\ref{kbounds}),
it follows from (\ref{a-bounds}) that
$|a(\btheta_0)|\leq 1$ when $n$ is sufficiently large.
It is routine to check that for any complex 
number $z$ with $|z|\leq 1$,
\[ \exp(-|z|) \leq \Re(\exp(z)) \leq \exp(|z|).\]
By~\eqref{a-bounds} we can apply this with $z = a(\btheta_0)$
to find that 
\begin{equation}
\int_{\R} G(\btheta) \, d\btheta = 
       \exp\(O(\lambda \, 2^n \Delta^3)\)\, 
                   \int_{\R} \exp \biggl( -\dfrac{1}{2}\, A 
                \sum_{\alpha\in\{ \pm 1\}^n} f_{\alpha}(\btheta)^2
                                       \biggr) \, d\btheta.
  \label{uptohere}
\end{equation}

Now we calculate that
\[
\sum_{\alpha \in \{ \pm 1\}^n}
   f_{\alpha}(\btheta)^2 = 2^n\sum_{S\in\S} \theta_S^2.
\]
Since this quantity is real and $\lambda 2^n\Delta^3 =o(1)$,
we have that
\[
  \int_\R G(\btheta)\,d\btheta =
  \(1 + O(\lambda 2^n\Delta^3)\)\prod_{S\in\S}
      \int_{-\Delta(2n)^{-|S|}}^{\Delta(2n)^{-|S|}}
        \exp\(-\dfrac12\lambda(1-\lambda)2^n\theta_S^2\)\,d\theta_S.
\]
Next we apply the well-known estimate
\[
  \int_{-x\sigma}^{x\sigma} e^{-u^2/(2\sigma^2)}\,du =
     \sigma\sqrt{2\pi}\,\(1 + o(e^{-x^2/2})\)
  \quad\text{for $x\to\infty$,}
\]
with $\sigma=\(\lambda(1-\lambda)2^n)^{-1/2}$
and $x=\Delta(2n)^{-|S|}\sigma^{-1}> \sqrt{32 nM}\to\infty$.
This gives
\[
\int_\R G(\btheta)\,d\btheta =
  \biggl( \frac{2\pi}{\lambda(1-\lambda)2^n}
           \biggr)^{\!\!M/2}
     \(1 + O(\lambda 2^n\Delta^3)+O(Me^{-16nM})\).
\]
The lemma follows on noting that the second error term is
subsumed by the first.
\end{proof}

\section{Proofs of Theorem \ref{main} and its corollaries}
\label{proofsection}

The theory we have developed over the preceding sections allows
us to complete the proofs of our main results.

\begin{proof}[Proof of Theorem~\ref{main}]
By~\eqref{Idef} we have that
\[
N(n,k,q)
=\frac{(1+\radius )^{2^n}}{(2\pi)^M\radius^{2^kq}}
\biggl(\,  \int_{\RRc }G(\btheta)\,d\btheta 
         + \int_{(\RRc)^c}G(\btheta)\,d\btheta
\biggr).
\]
First suppose that $\lambda\le\tfrac12$. 
Then the first integral is $2^Q\int_{\R}G$, where
$\int_{\R}G$ has been evaluated in Lemma~\ref{inthebox}, while
the second integral is bounded in absolute value by Lemma~\ref{outofbox}
and hence is covered by the error term of Lemma~\ref{inthebox}.
This completes the proof when $\lambda \le \tfrac12$, and the
result follows for $\lambda>\tfrac12$ by Remark~\ref{remark1}.
\end{proof}

Corollary~\ref{resilient} follows from Theorem~\ref{main} by
setting $\lambda=\tfrac12$.
Corollary~\ref{rightfirsttime} requires a little more effort.

\begin{proof}[Proof of Corollary~\ref{rightfirsttime}]
  We divide the interval of summation into five ranges.
  Define 
  \[ q_1=\lceil 2^{n-k-1} n^{-1} \rceil, \quad
 q_2=\lceil 2^{n-k-1}-2^{n/2-k}n\rceil, \quad
  q_3=2^{n-k}-q_2, \quad q_4=2^{n-k}-q_1.
  \]
  Also define
  \[
      W(\lambda)=W(\lambda,k,n) = 
     2^Q\(\pi A 2^{n+1}\)^{-M/2}
           \( \lambda^{\lambda }(1-\lambda)^{1-\lambda}\)^{- 2^n},
  \]
  which is the right side of~\eqref{theanswer} apart from the error term.

  We start with the range $q\in[q_2,q_3]$, for which
  $\lambda=\tfrac12+O(2^{-n/2}n)$.
  By Taylor expansion, we have for $x=O(2^{-n/2}n)$
  that
  \begin{equation}\label{Wexpansion}
     W\(\tfrac12+x\) 
       = W\(\tfrac12\)\exp\(-(2^{n+1}-2M)x^2+O(2^{-n}n^4)\).
  \end{equation}
  The error term in~\eqref{Wexpansion}
  is smaller than $2^{-n/2+3k}n^{3k+3/2}M^{3/2}$
  for any $\lambda$ in this range so, by Theorem~\ref{main},
  \[
      \sum_{q=q_2}^{q_3} N(n,k,q)
          = \(1 + O(2^{-n/2+3k}n^{3k+3/2}M^{3/2})\) W\(\tfrac12\)
              \sum_{q=q_2}^{q_3} h(q),
  \]
  where $h(q) = \exp\(-2^{-2n+2k+1}(2^n-M)(q-2^{n-k-1})^2\)$.
  By Euler-Maclaurin summation (see for example~\cite[p.\,36]{Wong}),
  \begin{align*}
      \sum_{q=q_2}^{q_3} h(q) &= O(e^{-n^2}) 
             + \(1 + O(2^{-n})\)\int_{q_2}^{q_3} h(q)\,dq\\
       &= \(1 + O(2^{-n})\)\pi^{1/2}2^{n-k-1/2}(2^n-M)^{-1/2}\\
       &= \(1 + O(2^{-n}M)\) \pi^{1/2}2^{n/2-k-1/2}.
  \end{align*}
  This proves that $\sum_{q=q_2}^{q_3} N(n,k,q)$ is given
  by an expression of the same form
  as the right side of~\eqref{Nsum}.

  Next consider the range $q\in [q_1,q_2)$, which is the
  mirror image of
  the range $q\in (q_3,q_4]$.   Then
  \begin{align*}
     \frac{d\log W(\lambda)}{d\lambda} &= 
       a_1(\lambda)2^n
        + a_2(\lambda)\(\lambda(1{-}\lambda)2^n-M\),
        \mbox{~~where}\\
     a_1(\lambda)&=\log(\lambda^{-1}-1)+\lambda-\dfrac12
     \mbox{~~and~~}
     a_2(\lambda)=\frac{1-2\lambda}{2\lambda(1{-}\lambda)}.
  \end{align*}
  For $\tfrac{1}{2n} \le\lambda\le 1-\tfrac{1}{2n}$ we
  find that $\lambda(1{-}\lambda)2^n>M$, while
  $a_1(\lambda)$ and $a_2(\lambda)$ have the same sign as
  $\tfrac12-\lambda$.  Therefore $W(\lambda)$ is unimodal
  in this range.  Since $\eta(n,k,q)=o(1)$, we have
  \begin{align*}
     \sum_{q=q_1}^{q_2-1} N(n,k,q)
        &= O(2^n) W\(\tfrac12-2^{-n/2}n\)\\
  &= O(2^n)\, \exp\( - (2 - o(1)) n^2\) \, W\(\tfrac12\)\\
  &= e^{-O(n^2)}\, W\(\tfrac12\) 
  \end{align*}
  using~\eqref{Wexpansion}.  This shows that the sum
over $[q_1,q_2)$ is covered by the error term of the Corollary.
By Remark~\ref{remark1} the same conclusion holds for the summation
from $q_3+1$ to $q_4$.
  
  Finally consider the range $q\in [0,q_1)$, which is the mirror image
  of the range $q\in (q_4,2^{n-k}]$.   Here
we use the trivial bound
  \[
     \sum_{q=0}^{q_1-1} N(n,k,q)
      < \(2^n\)^{2^k q_1} = 2^{2^{n-1} + O(n2^k)}
  \]
  which also fits into the error term of the corollary.
By Remark~\ref{remark1}, the same conclusion holds for the summation
from $q_4+1$ to $2^{n-k}$, which completes the proof.
\end{proof}

\section{More on the case $k=1$}
\label{k=1section}

In the case of $k=1$, which corresponds to the ``balanced colourings''
enumerated by Palmer, Read and Robinson~\cite{Palmer}, it is possible
to fill in the range of very small or very large values of~$q$ excluded
by~\eqref{kbounds}.

\begin{lemma}\label{smallq}
If\/ $0\le q=o(2^{n/2})$ then
\[
 (2q)!\, N(n,1,q) = \binom{2q}{q}^{\!n} \( 1+ O(q^2/2^n) \).
\]
\end{lemma}

\begin{proof} Generate a
$2q\times n$ matrix by a random process:
for each column independently, randomly insert 0 in $q$ rows
and 1 in the other~$q$.
This matrix is one of those counted by $(2q)!\,N(n,1,q)$ provided all
the rows are different.  (Recall that $N(n,1,q)$ counts matrices up
to row order.)

The probability that a specified pair of rows are equal is
\[
   2^n\biggl(\binom{2q{-}2}{q}{\Bigm/}\binom{2q}{q}\biggr)^{\!n}
       = \biggl(\frac{q-1}{2q-1}\biggr)^{\!n} < \,2^{-n},
\]
so, by the Bonferroni inequality, the probability that no two
rows are equal is $1-O(q^2/2^n)$.  This completes the proof.
\end{proof}

\begin{theorem}\label{k=1}
   Uniformly for $0\le q\le 2^{n-1}$,
\[
    N(n,1,q) =  \binom{2^n}{2q}
      \left(\! \frac{\displaystyle\binom{2^{n-1}}{q}^{\!\!2}}
                  {\displaystyle\binom{2^n}{2q}}\right)^{\!\!\textstyle n}
      \(1 + o(n^52^{-n/5})\).
\]  
\end{theorem}

\begin{proof}
We begin by motivating the given formula.
Choose, uniformly at random, a set of $2q$ distinct
elements of $\{0,1\}^n$.
The event that exactly $q$ of these elements have 1
in some specified position has probability
$$
  \binom{2^{n-1}}{q}^{\!\!2}{\Bigm/}\binom{2^n}{2q}.
$$
Therefore, the theorem is stating that these $n$ events are very
close to being independent in some sense.

We can derive the theorem from Theorem~\ref{main} and
Lemma~\ref{smallq}.  First consider the case that
$2^{2n/5}n^{12/5}\le q\le 2^{n-2}$.  Then, by Stirling's
formula,
$$
   \Biggl(\binom{2^{n-1}}{q}^{\!\!2}{\Bigm/}\binom{2^n}{2q}
   \Biggr)^{\mkern-5mu n}
    = \(\pi A 2^{n-1}\)^{-n/2}
      \(1 + O(n/q)\)
$$
and
\[ \binom{2^n}{2q} = (\pi A 2^{n+1})^{-1/2}\, 
   \(\lambda^\lambda (1-\lambda)^{1-\lambda})^{-2^n}\, (1+ O(1/q))\]
and the theorem follows from Theorem~\ref{main}.

In the case that $0\le q\le 2^{2n/5}n^{12/5}$, we calculate
that
$$
   \binom{2^{n-1}}{q}^{\!\!2}{\Bigm/}\binom{2^n}{2q} 
    = \binom{2q}{q}\binom{2^n-2q}{2^{n-1}-q}{\Bigm/}\binom{2^n}{2^{n-1}}
    = \binom{2q}{q}2^{-2q}\(1 + O(q/2^n)\),
$$
so the theorem follows from Lemma~\ref{smallq}.

Finally, for $2^{n-2}\le q\le 2^{n-1}$, take the complement
as in Remark~\ref{remark1}, noting that the binomial coefficients
in the statement of the theorem are symmetric around $q=2^{n-2}$.
\end{proof}

\section{Final remarks}\label{finalsection}

As mentioned in the Introduction, Denisov in~\cite{Denisov2}
incorrectly repudiated the result from~\cite{Denisov1} that we
quoted as Theorem~\ref{denisov1}.
Denisov's mistake was due to the incorrect computation of the
matrix inverse $A^{-1}$ on page~95 of~\cite{Denisov2}.
In fact the $I,J$ element of $A^{-1}$ is
$(-1)^{|J|-|I|}2^{|I|}$ for $I\subseteq J$ and $0$ otherwise.
Correcting the mistake shows that the critical
value $\bar z^T Q^{-1}\bar z$ on page~97 equals $2^{2k-n+2}$
and not the value stated.  Except for this error,
Denisov would have extended Theorem~\ref{denisov1} to
$k=o(n^{1/2})$ and in fact would have matched
Theorem~\ref{main} (with a different
vanishing error term) for $k=o(n^{1/2})$ and
\[ \abs{q-2^{n-k-1}}< \rho \, 2^{n/2-k}\,n^{1/2}\] 
for any positive constant $\rho < \sqrt{\dfrac{\log 2}{2}}$.
Note that our coverage of both $k$ and $q$ is considerably
wider than that.

\medskip

Finally we mention a
connection between correlation-immune boolean function and
Hadamard matrices. 
Recall that a {\it Hadamard matrix of order} $n$ is an $n\times n$
matrix over $\pm 1$ whose columns are pairwise orthogonal.
Such matrices are known to exist for $n=1$, $n=2$, and for infinitely
many other $n$.
If $n>2$ then $n\equiv 0 \mod 4$ is a necessary condition for the
existence of a Hadamard matrix of order $n$.  It is a long-standing
open problem to show that this necessary condition is also sufficient.
Let $H_n$ be the number of Hadamard matrices of order $n$.  By a simple
normalization, it can be seen that $H_n$ equals $2^n$ times the number
of Hadamard matrices whose leftmost column equals all $+1$'s.
If such a column is removed, and each $-1$ changed to $0$,
there remains an $n\times(n{-}1)$ matrix 
of the sort counted (up to row permutation) by $N(n{-}1,2,n/4)$.
Hence, $H_n = 2^n n!\, N(n{-}1,2,n/4)$ for $n > 2$.
This connection
raises the possibility of proving the Hadamard conjecture by
asymptotic methods.  Unfortunately, the coverage of Theorem~\ref{main}
is inadequate for that purpose.



\begin{thebibliography}{99}

\bibitem{Bach}
E.~Bach,
Improved asymptotic formulas for counting correlation-immune
Boolean functions,
Technical Report 1616, Computer Sciences Dept.,
University of Wisconsin, 2007.

\bibitem{Carlet-chapter}
C.~Carlet, 
Boolean functions for cryptography and error correcting codes, Preprint.  
To appear as a chapter of 
\emph{Boolean Functions: Theory, Algorithms and Applications} 
(Y.~Crama and P.~Hammer, eds.), Cambridge University Press.

\bibitem{Carlet}
C.~Carlet and A.~Gouget,
An upper bound on the number of $m$-resilient Boolean functions,
ASIACRYPT 2002, {\it Lecture Notes in Comput. Sci.}
{\bf 2501} (2002) 484--496.

\bibitem{CK}
C.~Carlet and A.~Klapper,
Upper bounds on the number of resilient functions and of bent functions, 
Springer-Verlag, Lecture Notes dedicated
to Philippe Delsarte (to appear). A shorter version has appeared
in the Proceedings of the 23rd Symposium on Information Theory in
the Benelux, Louvain-La-Neuve, Belgian, 2002.

\bibitem{CS}
T.\,W.~Cusick and P.~Stanica,
Cryptographic Boolean Functions and Applications,
Elsevier, 2009.

\bibitem{Denisov1}
O.\,V.~Denisov,
An asymptotic formula for the number of correlation-immune
of order~$q$ boolean functions, {\it Discrete Math. Appl.},
{\bf 2} (1992) 279--288.
originally published in {\it Diskretnaya Matematika} {\bf 3} (1990)
25--46 (in Russian).  Translated by A.V.~Kolchin.

\bibitem{Denisov2}
O.\,V.~Denisov,
A local limit theorem for the distribution of a part of the spectrum
of a random binary function,
{\it Discrete Math. Appl.}, {\bf 10} (2000) 87--101.
originally published in {\it Diskretnaya Matematika}, {\bf 12},1 (2000)
 (in Russian).  Translated by the author.

\bibitem{Hedayat}
A.\,S.~Heydayat, N.\,J.\,A.~Sloane, J.~Stufken,
Orthogonal arrays : theory and applications,
Springer-Verlag, 1999.

\bibitem{Maitra}
S.~Maitra and P.~Sarkar,
Enumeration of correlation immune boolean functions,
ACSIP'99, {\it Lecture Notes in Comput. Sci.},
{\bf 1587} (1999) 12--25.

\bibitem{mitchell}
C.~Mitchell, 
Enumerating Boolean functions of cryptographic significance, 
\emph{J. Cryptology}, {\bf 2} (1990), 155-170.

\bibitem{Palmer}
E.\,M.~Palmer, R.\,C.~Read and R.\,W.~Robinson,
Balancing the $n$-cube: a census of colorings,
{\it J. Algebraic Combin.}, {\bf 1} (1992) 257--273.

\bibitem{PLSK}
S.\,M.~Park, S.\,J.~Lee, S.\,H.~Sung and K.\,J.~Kim, 
Improving bounds for the number of correlation immune Boolean functions, 
\emph{Inform. Process. Lett.}, {\bf 61} (1997), 209--212.

\bibitem{Roy}
N.~Roy,
A brief outline of research on correlation immune functions,
ACISP 2002, {\it Lecture Notes in Comput. Sci.},
{\bf 2384} (2002) 379--394.

\bibitem{Sarkar}
P.~Sarkar,
A note on the spectral characterization of correlation immune
Boolean functions,
{\it Inform. Process. Lett.}, {\bf 74} (2000) 191--195.

\bibitem{Schneider}
M.~Schneider, A note on the construction and upper bounds of
correlation-immune functions, {\it Lecture Notes in Comput. Sci.},
{\bf 1355} (1997) 295--306.

\bibitem{tarannikov}
Y.~Tarannikov,
On the structure and numbers of higher order correlation-immune functions, 
in Proceedings of IEEE International Symposium on Information Theory, 2000, 185. 

\bibitem{TK}
Y.~Tarannikov and D.~Kirienko,
Spectral analysis of high order correlation immune functions,
Proceedings of 2001 IEEE International Symposium on Information Theory, 2001, 69.

\bibitem{Wong}
R.~Wong, \textit{Asymptotic approximations of integrals},
Academic Press, Boston, 1989.

\bibitem{Xiao}
G-Z.~Xiao and J.\,L.~Massey,
A spectral characterization of correlation-immune
combining functions.
{\it IEEE Trans. Inform. Theory}, {\bf 34} (1988) 569--571.

\bibitem{YG}
Y.\,X.~Yang and B.~Guo, 
Further enumerating Boolean functions of cryptographic significance, 
\emph{J. Cryptology}, {\bf 8} (1995), 115--122.

\bibitem{Zhang}
J-Z.~Zhang, Z-S.~You and Z-L.~Li,
Enumeration of binary orthogonal arrays of strength~1,
{\it Discrete Math.}, {\bf 239} (2001) 191--198.

\end{thebibliography}
\end{document}